\documentclass[reqno]{amsart}
\usepackage{amssymb}
\usepackage{amscd}
\usepackage{mathrsfs,amsmath}
\usepackage{enumitem}
\usepackage{xcolor}

\newtheorem{theorem}{Theorem}[section]

\newtheorem{lemma}{Lemma}[section]

\newtheorem{proposition}{Proposition}[section]
\newtheorem{remark}{Remark}[section]

\numberwithin{equation}{section}

\usepackage{xcolor}

\usepackage{empheq}

\newcommand*\interior[1]{\mathring{#1}}

\title[Existence result for the coupled problem]{Existence Result for a Model Coupling a Quasi-Linear Parabolic Equation and a Linear Hyperbolic System}

\subjclass[2010]{35A01, 35K20, 35L20, 35L53.}

\keywords{Hyperbolic system; Stokes equation; Fluid-Structure Interaction; Inverse energy estimate; Fixed point theory.}

\author[D. AIT-AKLI]{\bfseries Djamal AIT-AKLI$^{1,*}$\\\today}

\address{$^1$L2CSP laboratory, Mouloud Mammeri University of Tizi-Ouzou, Tizi-Ouzou, 15000, Algeria.}
\address{$^*$Corresponding author: \textnormal{djamel.aitakli@ummto.dz}}

\begin{document}


\begin{abstract}
We prove globally-in-time existence of solution for a problem coupling the linear Lam\'e system and the quasi-linear Stokes equation. A solution of this global coupled problem is viewed as the fixed point of some non-linear operator $T$. We construct, using a regularization procedure, a sequence $(T^\epsilon)_\epsilon$ of auxiliary approximating compact operators. Then we establish, using a combination of Banach and Schaeffer fixed point theorems, the existence of fixed points to every operator $T^\epsilon$. Finally we prove that these fixed points converge to the fixed point of $T$. 
\end{abstract}

\maketitle

\section{Introduction}    
We address the issue of existence of solution for the coupled system which reads:
\begin{empheq}[left={\empheqlbrace}]{alignat=2}
&\partial_t v_f-\nu {\rm div}a(t, \nabla v_f) -\nu\Delta v_f + \nabla \pi = F ~\text{in}~\Omega_f^T,\label{sys1eq1}\\
&\mathop{\rm div}v_f = 0 ~\text{in}~\Omega_f^T,\quad v_f = 0\quad \text{on}\quad(0,T)\times(\partial\Omega_f - \Sigma),\label{sys1eq2}\\
& v_f(0,.) = v_f^0\quad \text{in}\quad\Omega_f,\label{sys1eq3}\\
&\int_0^t v_f(r){\rm d}r = u_s(t) , \quad \left( S(v_f,\pi) + a(t, \nabla v_f)\right)\cdot\overrightarrow{n} =  \sigma(u_s)\cdot \overrightarrow{n}\quad\text{on}\quad \Sigma^T,\label{sys1eq4}\\
&\partial_{tt} u_s -\mathop{\rm div}\sigma(u_s)= 0  ~\text{in}~ \Omega_s^T,\label{sys1eq5}\\
&u_s = 0  \quad\text{on}~ (0,T)\times(\partial\Omega_s - \Sigma),\label{sys1eq6}\\
&u_s(0,.) = 0, \quad \partial_t u_s(0,.) = 0  ~\text{in}~\Omega_s.\label{sys1eq7}
\end{empheq}
Let us describe the components and the notations in system (\ref{sys1eq1})-(\ref{sys1eq7}). We let $T>0$ to be a positive real number. We have denoted $ \omega^T := (0,T)\times\omega$. The function $u_s$ denotes the displacement of the solid structure, the function $v_f$ denotes the velocity of the fluid  and $\pi$ denotes its pressure. System  (\ref{sys1eq1})-(\ref{sys1eq7}) is formed out of:
\begin{itemize}
\item The parabolic Stokes equation with a quasi-linear diffusion term that describes the motion of a fluid inside a fluid domain $\Omega_f\subset\mathbb{R}^2$, cf. equations (\ref{sys1eq1})-(\ref{sys1eq3}).
\item The coupling condition, cf. equation (\ref{sys1eq4}). 
\item The second order linear hyperbolic Lam\'e system that describes the deformation, due to interaction with the fluid, of a structure occupying a solid domain $\Omega_s\subset\mathbb{R}^2$, cf. equations (\ref{sys1eq5})-(\ref{sys1eq7}).
\end{itemize}

\noindent We assume that $\Omega := \Omega_s\cup\Omega_f(t)\subset\mathbb{R}^2$, for all $t\in (0,T)$, that is the global domain $\Omega$ doesn't vary in time during the interaction. Moreover, the boundaries $\partial\Omega_f$ resp. $\partial\Omega_s$ of the fluid resp. the structure domains are both assumed to meet the minimal $C^2-$regularity. 

The solid domain $\Omega_s$ and the fluid domain $\Omega_f(t)$ share a common part of their respective boundaries, this contact interface is denoted $\Sigma$ i.e. $\partial\Omega_s\cap\partial\Omega_f(t) = \Sigma(t)$, the variation in time of $\Sigma$ isn't taken into account, cf. the conclusion at the end of this paper, we assume furthermore that $\Sigma$ is connected. The coupling condition, relating the fluid and solid problems, is prescribed on $(0,T)\times\Sigma$. It consists of imposing, on the contact interface $\Sigma$, the two equalities:
\begin{equation}\label{eq1}
  \begin{array}{ccl}
  &\int_0^t v_f(r,x){\rm d}r = u_s(t,x)  \\
    &\left( S(v_f,\pi) + a(t, \nabla v_f)\right)\cdot\overrightarrow{n} =  \sigma(u_s)\cdot \overrightarrow{n},\quad (t,x)\in\Sigma^T,
  \end{array}
\end{equation}
where $\overrightarrow{n}$ denotes the exterior unit normal defined at each point of $\Sigma$. The term $S(v_f,\pi)$ denotes the Cauchy stress tensor:
\begin{equation}\label{eq2}
S(v_f, \pi) :=   - \pi \mathbb{I} + 2\nu \varepsilon(v_f),
\end{equation}
and the function $a: [0,T]\times\mathbb{R}^2\to \mathbb{R}^2$ such that $a := (a_1, a_2)$ satisfy (\ref{monotone}) and other assumptions that will be precised in section \ref{sec3}. On the other hand, $\sigma(u_s)\cdot\overrightarrow{n}$ denotes the normal component of the stress tensor:
\begin{equation}\label{eq3}
\sigma(u) = 2\mu\varepsilon(u) + \lambda \textrm{ Tr}\varepsilon(u) \textrm{ Id} ,\quad\text{with}\quad  \varepsilon(u) = \frac{1}{2}\left( \nabla u + \nabla^t u  \right).
\end{equation}
Unlike what is usually done regarding the coupling condition, we imposed the equality of the displacements functions, $u_s$ and $ u_f$, on the contact interface instead of the equality of the velocities.
  
 The condition (\ref{eq1}) is completed by an homogeneous Dirichlet data (\ref{sys1eq2}), (\ref{sys1eq6}) on the remaining part of the boundaries of the domain and by initial time conditions (\ref{sys1eq3}) and (\ref{sys1eq7}). The fluid is assumed to be divergence-free.
  We emphasize that, except for the fluid initial condition, (\ref{sys1eq3}), the restriction to homogeneous data for both Stokes and Lam\'e problems are adopted only for the sake of simplicity of presentation, one can refer to \cite{1} and the references therein for the case where non homogeneous data are considered but with restrictive assumption on the  geometry of the domain, namely the flatness of the contact interface.
  
   It is worth noting that this kind of coupled system of PDEs has proven efficient for modeling various Fluid-Structure interactions (FSI) phenomenon.  The main result of the present work is stated in the following theorem:
\begin{theorem}\label{main-The}
$\forall T>0$, $\forall F\in L^2(0,T; H^{-1}(\Omega_f))$ and $\forall v_f^0\in L^2(\Omega_f)$, the coupled problem (\ref{sys1eq1})-(\ref{sys1eq7}) admits at least one solution:
\begin{equation}\label{RegularityOfSolution}
\begin{aligned}
\left(v_f, \pi, u_s\right)\in &L^2\left(0,T,  H^1(\Omega_f)\right)\cap H^1(0,T; H^{-1}(\Omega_f)\times L^2\left( 0,T; L^2(\Omega_f)/\mathbb{R}\right)\\&\times\left( H^1((0,T)\times\Omega_s)\cap H^2(0,T; H^{-1}(\Omega_s))\right).
\end{aligned}
\end{equation}
\end{theorem}
\noindent The focus is on establishing existence of globally-in-time solution to the coupled problem (\ref{sys1eq1})-(\ref{sys1eq7}) which is is analogous to  problem {\cite[Problem 2.10, p.556]{1}} 
The interest in the regularity of the solutions is secondary in our current considerations. 
 
\noindent The novelty in this work is about two main things. The first one consists in proving the globally-in-time existence of at least one solution to the coupled problem while dealing with the non-linearity in the Stokes equation. The second new thing consists in considering domains with arbitrary geometry, that is no particular restriction is assumed on the contact surface $ \Sigma $, regarding namely its geometric flatness, this extends the result proved by the authors in \cite{1}. To achieve that, we view the solution of (\ref{sys1eq1})-(\ref{sys1eq7}) as a fixed point of some non-linear operator $T$ and we use a regularization method in order to apply fixed point theory.
\subsection*{Organization of the paper}
	In the second section we establish a well-posedness result for the Dirichlet problem associated to the Lam\'e operator along with an inverse estimate. In the third section we derive an energy estimate for the quasi-linear Stokes system and we introduce the operator $T: \mathcal{X}\rightarrow\mathcal{X}$ whose fixed point is a solution of problem (\ref{sys1eq1})-(\ref{sys1eq7}), the space $\mathcal{X}$ is given by (\ref{X}). Next we apply the regularization method to construct a sequence of auxiliary compact operators $(T^\epsilon)_\epsilon$ such that $T =\underset{\epsilon\to 0}{\lim} T^\epsilon$, then we establish the boundedness and compactness of $T^\epsilon$ using the preceding estimates. In the fourth section, we prove the existence of a fixed point $u^\epsilon$ to $T^\epsilon$ by combining Banach and Schaeffer fixed point theorems. Finally we conclude by showing that the fixed points $u^\epsilon$ converge to a fixed point $u^0$ of $T$.

\section{Inverse estimate for the Lam\'e system.}
Throughout this section, we let $\Omega_s\subset\mathbb{R}^2$ to be a bounded planar domain  with boundary $\partial\Omega_s\in C^2$. We consider the auxiliary Dirichlet problem (\ref{sys2}) associated to the time dependent second order Lam\'e operator $\mathcal{H}$ given by:
\begin{equation}\label{H}
\mathcal{H}u := \partial_{tt}u -{\rm div}\sigma(u).
\end{equation}
We prescribe a non-homogeneous Dirichlet condition on $(0,T)\times\Sigma$ and a homogeneous Dirichlet condition on the remaining part of the boundary: 
\begin{equation}\label{sys2}
\left\{
    \begin{array}{ll}
\begin{aligned}
& \partial_{tt}u_s  - {\rm div}\sigma(u_s) = 0\quad\text{in}~(0,T)\times\Omega_s,\\
& u_s = u_s^d \quad\text{on}~ (0,T)\times\Sigma,\\
& u_s = 0  \quad\text{on}~ (0,T)\times \left(\partial\Omega_s -\Sigma \right),
\end{aligned}
   \end{array}
\right. 
\end{equation}
where $T>0$. Moreover, initial-time conditions are prescribed: 
\begin{equation}\label{eq4}
u_s(0) = 0,\quad \partial_tu_s(0) = 0.
\end{equation} 
The Dirichlet data $u_s^d$ in (\ref{sys2}) is assumed to be compatible with (\ref{eq4}).
Let us define the space: $$\mathcal{D}_s:= \lbrace v\in C^\infty(\overline{\Omega_s}),~ {\rm supp}v\cap(\partial\Omega_s - \interior{\Sigma})= \emptyset \rbrace,$$ where $\interior{\Sigma}$ denotes the topological interior of $\Sigma$. We set: 
\begin{equation}\label{U}
\mathcal{U}:= \overline{\mathcal{D}_s}^{H^1}
\end{equation}
 to be the completion of $\mathcal{D}_s$ with respect to the $H^1(\Omega_s)-$norm. Consider the space:
\begin{equation}\label{equivalence}
\lbrace \gamma_\Sigma(v):~ v\in\mathcal{U} \rbrace \equiv H_0^\frac{1}{2}(\Sigma),
\end{equation}
where $\gamma_\Sigma(v)$ denotes the trace, on the boundary $\Sigma$, of the function $v$. The space $H_0^\frac{1}{2}(\Sigma)$ is defined as the completion of $C_0^\infty(\Sigma)$ with respect to the $H^\frac{1}{2}(\Sigma)-$norm. Denote $\mathcal{X}$ to be the space: 
\begin{equation}\label{X}
\mathcal{X} := H^\frac{1}{2}\left(0,T; L^2(\Sigma)\right)\cap L^2\left(0, T; ~H_0^\frac{1}{2}(\Sigma)\right)\subset H^\frac{1}{2}((0,T)\times\Sigma),
\end{equation}
and denote $[\mathcal{X}]^\ast$ its topological dual space. Let us define the operator $T_1$ by:
\begin{equation}\label{T1}
\begin{aligned}
T_1 : \mathcal{X} &\rightarrow L^2\left(0, T; H^{-\frac{1}{2}}(\Sigma)\right)\\
                 u_s^d  &\mapsto g_s = T_1(u_s^d).
\end{aligned} 
\end{equation} 
where $H^{-\frac{1}{2}}(\Sigma)$ denotes the dual of $H_0^\frac{1}{2}(\Sigma)$. The operator $T_1$ associates to every Dirichlet data $u_s^d\in\mathcal{X}$ on the solid part of the contact interface, $(0,T)\times\Sigma$, the uniquely determined Neumann data $g_s := \sigma(u_s)\cdot\overrightarrow{n}$ corresponding to the solution $u_s$ of the Dirichlet problem (\ref{sys2})-(\ref{eq4}). In the sequel we will denote equally by $u_s$ the Dirichlet data $u_s^d$. The main result of this section is  given in the following proposition:

\begin{proposition}\label{prop1}
The operator $T_1$, given by (\ref{T1}), is well defined and bounded i.e. there exists a constant $C_s>0$ such that:
\begin{equation}\label{estimlame}
||\sigma(u_s)\cdot\overrightarrow{n}||_{L^2\left( 0,T; ~ H^{-\frac{1}{2}}(\Sigma)\right)} \leq C_s ||u_s||_{\mathcal{X}}
\end{equation}
for all $u_s\in \mathcal{X}$, where $\mathcal{X}$ is defined by (\ref{X}).
\end{proposition}
Before passing on to the proof of Proposition \ref{prop1}, we give a remark stating a lifting property within the context of the Bochner space $L^2(0,T; \mathcal{U})$ :
\begin{remark}\label{remark}
Consider the map $\gamma_\Sigma$ defined by:
\begin{equation*}
\begin{aligned}
\gamma_\Sigma:   L^2(0,T; \mathcal{U})&\rightarrow L^2\left(0,T; H_0^{\frac{1}{2}}(\Sigma)\right)\\
                 v  &\mapsto \gamma_\Sigma v,
\end{aligned} 
\end{equation*}
this map associates to every $v\in  L^2(0,T; \mathcal{U})$ its trace on $(0,T)\times\Sigma$. We claim that $\gamma_\Sigma$ is onto. Indeed, let $v\in L^2\left(0,T; H_0^{\frac{1}{2}}(\Sigma)\right)$. The function can be extended by zero, to the rest of the boundary, into a function in $L^2\left(0,T;H^{\frac{1}{2}}(\partial\Omega_s)\right)$ which we still denote by $v$. By using the lifting property, one can easily find a family of functions $ (0,T)\ni t\mapsto\tilde{v}(t,.)\in H^1(\Omega_s)$ such that $\gamma_{\partial\Omega_s}\tilde{v}(t,.) = v(t,.)$, $\forall t\in(0,T)$, and such that one also has: $||\tilde{v}(t,.)||_{H^1(\Omega_s)}\leq ||v(t,.) ||_{H^{\frac{1}{2}}(\Omega_s)}$, $\forall t\in(0,T)$. The existence of a $\tilde{v}$ satisfying such estimate can be established, for instance, within the Banach space $\lbrace w(t,.)\in H^1(\Omega_s):~ \int_{\Omega_s} w(t,x)\Delta_x\phi(x){\rm d}x =0, \forall\phi\in C^\infty_c(\Omega_s), \forall t\in (0,T) \rbrace$. Thus the function $\tilde{v}$ satisfy:
\begin{equation*}
||\tilde{v}||_{L^2(0,T; \mathcal{U})}\leq ||v||_{L^2\left( 0,T;H_0^{\frac{1}{2}}(\Sigma)\right)},
\end{equation*}
and moreover, by combining {\cite[Theorem 1, p.518]{2}} and {\cite[Theorem 2.1, p.731]{3}}, we have $\tilde{v}(t,.)\in\mathcal{U}$, $\forall t\in (0,T)$; we thus conclude the surjectivity of $\gamma_\Sigma$.
\end{remark}

Let us state a lemma about a useful existence and regularity result:

\begin{lemma}\label{lem22}
Consider the problem 
\begin{equation}\label{sys44}
\left\{
    \begin{array}{ll}
\begin{aligned}
& \partial_{tt}\phi  - {\rm div}\sigma(\phi) = F\quad\text{in}~(0,T)\times\Omega_s,\\
& \phi = 0  \quad\text{on}~ (0,T)\times \partial\Omega_s,\\
& \phi(0,.) = 0, \quad \partial_t\phi(0,.)=0,
\end{aligned}
   \end{array}
\right.
\end{equation}
we claim that for every $F\in L^2(0,T; H^{-1}(\Omega_s)) $,  problem (\ref{sys44}) admits a unique solution $\phi\in L^2(0,T; H^1(\Omega_s))$ and moreover $\sigma(\phi)\cdot\overrightarrow{n} \in L^2(0,T; H^{-\frac{1}{2}}(\partial\Omega_s))$.
\end{lemma}  

\begin{proof}
One can, by density, find a sequence $(F_n)_n$ of elements in $L^2(0,T; L^2(\Omega_s)) \subset L^1(0,T; L^2(\Omega_s))$ such that: 
\begin{equation}\label{convF}
||F - F_n||_{L^2(0,T; H^{-1}(\Omega_s))} \to 0.
\end{equation}
 Next one can use {\cite[Theorem 2.1, p.151]{10}} to show the existence of a unique solution $\phi_n\in L^2(0,T; H^1(\Omega_s))=: \Lambda$ to problem (\ref{sys44}) with $F_n$ as a right hand side instead of $F$. Starting from the weak formulation of  problem (\ref{sys44}), one writes:
\begin{equation}
<\partial_t \phi_n, \partial_t\psi>_{} + <\varepsilon(\phi_n), \varepsilon(\psi)>_{L^2(0,T; L^2(\Omega_s))} = <F_n,\psi>_{\Lambda^*, \Lambda}   
\end{equation}
for all $\psi \in L^2(0,T; H^1_0(\Omega))$. Choosing $\psi = \phi_n$ and using  (\ref{convF}) we infer easily that
\begin{equation}
||\partial_t\phi_n||_{L^2(0,T; L^2(\Omega_s))} + ||\varepsilon(\phi_n)||_{L^2(0,T; L^2(\Omega_s))} < C
\end{equation}
for some $C>0$, and deduce the existence of a unique solution $\lim_n \phi_n  := \phi\in \Lambda $ to problem (\ref{sys44}).
\end{proof}
Now we pass into the proof of Proposition \ref{prop1}:

\begin{proof} 
 Define the subspace 
\begin{equation}\label{A}
\mathcal{A} :=
\begin{cases}
\begin{rcases}  
v\in H^1((0,T)\times\Omega_s)\cap H^2(0,T; H^{-1}(\Omega_s)):&  \\
\int_0^T\int_{\Omega_s}v\mathcal{H}(\phi){\rm d}x{\rm d}t = 0,\forall\phi\in C^\infty_c((0,T)\times\Omega_s ),&  \\
||v(0)||_{L^2(\Omega_s)} = ||\partial_tv(0)||_{H^{-1}(\Omega_s)} = 0,&
\end{rcases}
\end{cases} 
\end{equation}
where $\mathcal{H}$ is defined in (\ref{H}). Recall that: $ H^1(0,T; L^2(\Omega_s))\cap L^2(0,T; H^1(\Omega_s)) =  H^1((0,T)\times\Omega_s).$
One sees that $\mathcal{A}$ is a Banach space when endowed with the norm: 
$$||v||_{\mathcal{A}} := ||v||_{H^1((0,T)\times\Omega_s)} + ||\partial_{tt}v||_{L^2(0,T; H^{-1}(\Omega_s))},$$moreover $\mathcal{A}$ is reflexive. The idea of the proof consists at writing $T_1$ as a composition  $T_1 = N\circ \gamma_0^{-1}$ of two linear operators and then establishing their boundedness. The rest of the proof is divided into two main steps:\\
{\bf Step 1: first inverse estimate}.
 Consider the trace operator:
\begin{equation*}
\begin{aligned}
\gamma_0: \mathcal{A}&\rightarrow H^{\frac{1}{2}}((0,T)\times\partial\Omega_s)\\
                 u_s  &\mapsto \gamma_0 u_s,
\end{aligned} 
\end{equation*}
this operator associates to every function $u_s\in \mathcal{A}$ its trace on $(0,T)\times\partial\Omega_s$, the space $\mathcal{A}$ is given by (\ref{A}). The linear operator $\gamma_0$ is clearly one-to-one, we claim that it is also onto. Indeed, let $u_s\in H^{\frac{1}{2}}((0,T)\times\partial\Omega_s)$, we are going to show the existence of $U_s \in \mathcal{A}$ such that $ u_s = \gamma_0 U_s$.
To do this, it suffices to prove that the problem:
\begin{equation}\label{sys3}
\left\{
    \begin{array}{ll}
\begin{aligned}
& \partial_{tt}U_s  - {\rm div}\sigma(U_s) = 0\quad\text{in}~(0,T)\times\Omega_s,\\
& U_s = u_s  \quad\text{on}~ (0,T)\times \partial\Omega_s,\\
& U_s(0,.) = 0, \quad \partial_tU_s(0,.)=0,
\end{aligned}
   \end{array}
\right.
\end{equation} 
admits a solution $U_s \in\mathcal{A}$. Let $u_s^n \in C^\infty((0,T)\times\partial\Omega_s)$ be such that: 
 \begin{equation}\label{convusn}
||u_s^n - u_s||_{H^{\frac{1}{2}}((0,T)\times\partial\Omega_s)}\to 0;
\end{equation}
such functions do exist by a density argument. Denote $U_s^n$ to be the solution of problem (\ref{sys3}) corresponding to $u_s^n$ as a Dirichlet data, according to {\cite[Theorem 2.1, p.151]{10}}, this problem admits a unique solution $U_s^n\in\mathcal{A}$ . The first equation of (\ref{sys3}) yields us:
\begin{equation}
< U^n_s, \partial_{tt}\phi>  - <U_s^n, {\rm div}\sigma(\phi)> = < U^n_s, \sigma(\phi)\cdot\overrightarrow{n}>,
\end{equation} 
 $\forall\phi\in L^2(0,T; H^1(\Omega_s))$ solution of problem (\ref{sys44}). Since, by Lemma \ref{lem22}, $\sigma(\phi)\cdot\overrightarrow{n} \in L^2(0,T; H^{-\frac{1}{2}}(\partial\Omega_s))$ then using the convergence (\ref{convusn})  we have: 
\begin{equation}\label{convv2}
< U^n_s, \mathcal{H}(\phi)>   = < U^n_s, \sigma(\phi)\cdot\overrightarrow{n}> \to   < u_s, \sigma(\phi)\cdot\overrightarrow{n}>,
\end{equation} 
for every $\phi\in L^2(0,T; H^1(\Omega_s))$ solution of problem  (\ref{sys44}) i.e. 
\begin{equation}\label{218}
< U^n_s, F>\quad {\rm converges} \quad \forall F \in L^2(0,T; H^{-1}(\Omega_s)). 
\end{equation}
Given that $L^2(0,T; H^{-1}(\Omega_s))$ is reflexive, then (\ref{218}) implies: 
\begin{equation}\label{uniformb1}
\exists C>0, \quad ||U_s^n||_{L^2(0,T; H^1(\Omega_s))}\leq C,\quad\forall n,
\end{equation}
and thus $(U_s^n)_n$ converges weakly to $\lim_{n\to\infty} U_s^n =: U_s \in L^2(0,T; H^1(\Omega_s))$. Furthermore, we have: 
\begin{equation}
<\partial_{tt} U^n_s, \phi> = <U_s^n, \partial_{tt}\phi>;
\end{equation}
 letting $n\to\infty$, we obtain
\begin{equation}
<\partial_{tt} U_s, \phi> = <U_s, \partial_{tt}\phi> ,\quad <\partial_t U_s, \partial_t\phi> = <U_s, \partial_{tt}\phi>,
\end{equation}
for all $\phi\in\mathcal{A}$, which implies that $U_s\in\mathcal{A}$.

Following the continuity argument stated in {\cite[Theorem 2.3, p.153]{10}}, we see that $||\partial_t U_s^n||_{L^2(0,T; H^{-1}(\Omega_s))}$ is uniformly bounded. Consequently, by combining this last fact and (\ref{uniformb1}) we deduce, by invoking the Aubin-Lions lemma and the continuity of the trace operator, that $||u_s^n - U_s ||_{L^2(0,T; L^2(\Sigma))} \leq||U_s^n - U_s||_{L^2(0,T; H^{\frac{1}{2}}(\Omega_s))}\to 0$,  this show that $\gamma_0 (\lim U_s^n) = \gamma_0(U_s) = u_s$, which concludes the surjectivity of $\gamma_0$. We thus infer that the operator $\gamma_0$ is an isomorphism. Since the domain and codomain  of the operator $\gamma_0$ are Banach spaces and since $\gamma_0$ is bounded, then by applying the Banach isomorphism theorem, we deduce that the inverse operator $\gamma_0^{-1}$ is bounded i.e. $\exists C_{\gamma_0^{-1}} >0$ such that: 
\begin{equation}\label{Ineq1}
||u_s ||_{\mathcal{A}}\leq C_{\gamma_0^{-1}} ||u_s  ||_{H^{\frac{1}{2}}((0,T)\times\Sigma)},
\end{equation}
for all $u_s \in\mathcal{X}\subset H^{\frac{1}{2}}((0,T)\times\Sigma)$.

\noindent{\bf Step 2: second inverse estimate}. Consider the following operator:
\begin{equation*}
\begin{aligned}
N : \mathcal{A}&\rightarrow L^2(0,T; H^{-\frac{1}{2}}(\Sigma))\\
                 u_s  &\mapsto N(u_s) = g_s,
\end{aligned} 
\end{equation*}
where $\mathcal{A}$ is defined by (\ref{A}). The operator $N$ associates to every displacement $u_s \in\mathcal{A}$ the corresponding Neumann data, $g_s:= \sigma(u_s)\cdot\overrightarrow{n}$, on the boundary $(0,T)\times\Sigma$. We claim that $N(\mathcal{A}) \subset L^2(0,T; H^{-\frac{1}{2}}(\Sigma)) $. Indeed, using the density of smooth functions in the space $\left(\mathcal{X}, ||~||_{H^\frac{1}{2}((0,T)\times\Sigma)}\right)$, then with the aid of estimate (\ref{Ineq1}) we can construct a sequence $(u^n_s)_n$ of elements in $C^\infty((0,T)\times\Omega)\cap\mathcal{A}$ such that:
\begin{equation}\label{convun}
||u^n_s - u_s||_{\mathcal{A}}\rightarrow 0 ~{\rm as}~ n\rightarrow \infty.
\end{equation}
Let $\phi\in L^2(0,T; \mathcal{U})$, where $\mathcal{U}$ is given by (\ref{U}). We integrate by part the first equation in (\ref{sys2}) against the test function $\phi(t,.)\in\mathcal{U}$ to obtain:
 \begin{equation}\label{formula1}
 \begin{aligned}
&\int_{0}^T<\partial_{tt} u_s^n, \phi>_{[1]^\ast, [1]}{\rm d}t +   \int_0^T\int_{\Omega_s} \varepsilon(u_s^n) \varepsilon(\phi)~{\rm d}x{\rm d}t\\& = \int_0^T<g_s^n , \phi  >_{[-\frac{1}{2}], [\frac{1}{2}], \Sigma}{\rm d}t,
\end{aligned}
\end{equation} 
for every $\phi \in  L^2(0,T; \mathcal{U})$, where $g^n_s := \sigma(u_s^n)\cdot\overrightarrow{n}\in C^\infty((0,T)\times\Sigma)$, the bracket $<., .  >_{[\beta]^\ast, \beta, E}$ denotes the duality pairing between $[H^{\beta}(E)]^\ast$ and $H^{\beta}(E)$. The tensor $\varepsilon$ is defined by (\ref{eq3}). Using (\ref{convun}), the expression (\ref{formula1}), the claim stated in Remark \ref{remark} and (\ref{equivalence}), we deduce that $\forall\phi \in L^2\left(0,T; H_0^\frac{1}{2}(\Sigma)\right)$:  
\begin{equation}\label{property1}
\begin{aligned}
&\left(|\int_0^T<g_s^n , \phi  >_{-\frac{1}{2}, \frac{1}{2}, \Sigma}{\rm d}t|\right)_n \text{ is a Cauchy sequence},
\end{aligned}
\end{equation}
the completeness of $\mathbb{R}$ yields $\sup_n|\int_0^T<g_s^n , \phi  >_{-\frac{1}{2}, \frac{1}{2}, \Sigma}{\rm d}t|<\infty$ for every $\phi\in L^2\left(0,T; H_0^\frac{1}{2}(\Sigma)\right)$. Using the uniform boundedness principle we infer that 
\begin{equation}\label{estim25}
\sup_n ||g_s^n||_{\mathcal{L}\left( L^2\left(0,T; H_0^\frac{1}{2}(\Sigma)\right), \mathbb{R}\right)} < \infty,
\end{equation}
 given the completeness and the separability of $L^2\left(0,T; H_0^\frac{1}{2}(\Sigma)\right)$, we infer using the fundamental theorem of weak$^*$ convergence and estimate (\ref{estim25}) that the sequence $(g^{\alpha(n)}_s)_n$ converges weakly$^*$ to some $g_s\in L^2\left(0,T; H^{-\frac{1}{2}}(\Sigma)\right)$, for some subsequence $(\alpha(n))_n$. Actually one can easily remark, using (\ref{property1}), that the whole sequence converges to $g_s$. We infer that the operator $N$ is well defined. The reader should notice that we have only proved: $N(\mathcal{A}) \subset L^2\left(0,T; H^{-\frac{1}{2}}(\Sigma)\right)$.

 Moreover, the operator $N$ is bounded. Indeed, knowing (\ref{equivalence}), we easily infer from the above arguments that the operator $N$ sends every weakly convergent sequence in $\mathcal{A}$ into a weakly$^*$ convergent sequence in $L^2\left(0,T; H^{-\frac{1}{2}}(\Sigma)\right)$. But given the reflexivness of the space $L^2\left(0,T; H_0^\frac{1}{2}(\Sigma)\right)$, the weak$^*$ convergence and the weak convergence agree. This shows that $N$ is sequentially weakly continuous. Since $N$ is linear, we deduce that it is bounded i.e. $\exists C_N>0$ such that:
\begin{equation}\label{Ineq2}
||g_s ||_{L^2\left(0,T; H^{-\frac{1}{2}}(\Sigma)\right)}\leq C_N ||u_s||_{\mathcal{A}},
\end{equation}
for every $u_s\in \mathcal{A}$.

Finally, by combining (\ref{Ineq1}) and (\ref{Ineq2}), we infer that  $\exists C_s>0$ such that:
\begin{equation}\label{estimff}
||\sigma(u_s)\cdot\overrightarrow{n}||_{L^2\left(0,T; H^{-\frac{1}{2}}(\Sigma)\right)} \leq C_s ||u_s||_{H^{\frac{1}{2}}((0,T)\times\partial\Omega_s)},
\end{equation}
for every $u_s\in H^{\frac{1}{2}}((0,T)\times\Sigma)$. Estimate (\ref{estimff}) holds for every $u_s\in \mathcal{X}\subset H^{\frac{1}{2}}((0,T)\times\Sigma) $, thus we conclude immediately estimate (\ref{estimlame}).
\end{proof}

\section{Estimates for the quasi-linear Stokes problem}\label{sec3}
\subsection*{Problem setting for the fluid part}
Assume $\Omega_f\subset\mathbb{R}^2$ to be a sufficiently smooth domain, say with boundary $\partial\Omega_f = \Sigma \cup( \partial\Omega_f-\Sigma)\in C^2$.
We consider the non-stationary Stokes operator with a quasi-linear diffusion term appearing in the first equation of (\ref{sys4}). This operator is endowed, cf. system (\ref{sys4}), with mixed boundary conditions. We prescribe a non-homogeneous Neumann condition on the contact interface, $(0,T)\times\Sigma$, via the Cauchy stress tensor, and  prescribe a homogeneous Dirichlet condition on the remaining part of the boundary. The fluid is assumed to be divergence-free. Let $F\in L^2(0,T ; H^{-1}(\Omega_f))$, $g_f\in L^2(0,T; H^{-\frac{1}{2}}(\Sigma)) $ and $v_f^0\in L^2(\Omega_f)$. The fluid part of the problem reads:

\begin{equation}\label{sys4}
\left\{
    \begin{array}{ll}
\begin{aligned}
&\partial_t v_f - \nu {\rm div~}a(t, \nabla v_f) - \nu \Delta v_f  +\nabla\pi  = F\quad\text{in}~(0,T)\times\Omega_f,\\
&{\rm div}v_f = 0  \quad\text{in}~ (0,T)\times\Omega_f,\\
&\left( a(t, \nabla v_f) + S(v_f,\pi)\right)\cdot\overrightarrow{n} = g_f \quad\text{on}~ (0,T)\times\Sigma,\\
&   v_f =  0 \quad\text{on}~(0,T)\times(\partial\Omega_f - \Sigma),\\
& v_f(0,.) = v_f^0 \quad\text{in}~ \Omega_f,
\end{aligned}
   \end{array}
\right.
\end{equation}
where $v_f$ is the unknown fluid velocity vector, $\pi$ denotes the unknown pressure and $S(v_f, \pi)$ 
denotes the Cauchy stress tensor given by (\ref{eq2}). For the viscosity, we assume for simplicity that $\nu =1$. 
 Let the vector function $a := (a_1, a_2)$ be such that the functions $a_j: [0, T]\times\mathbb{R}^2\rightarrow\mathbb{R}$, with $j=1,2 $, satisfy the assumptions stated in {\cite[Example 6.A, p.139]{5}} in the two dimensional case. One easily sees that these assumptions imply the hypothesis of {\cite[Proposition 5.1, p.129]{5}}. Actually we assume a stronger condition than {\cite[Condition 6.6.c, p.139]{5}}, that is: $\exists c_m>0$ such that:
\begin{equation}\label{monotone}
( a(t, \xi) - a(t, \eta))(\xi - \eta) \geq c_m|\xi - \eta|_2^2, \quad\forall\xi, \eta\in\mathbb{R}^2, \forall t\in[0,T].
\end{equation}
 Consider the space: 
\begin{equation}\label{Df}
\mathcal{D}_f := \lbrace v\in C^{\infty}(\overline{\Omega_f}),\quad {\rm div}~v =0, \quad \text{supp}~v\cap (\partial\Omega_f -\Sigma) = \emptyset\rbrace.
\end{equation}
Denote $\mathcal{V}$ to be the closure of $\mathcal{D}_f$ with respect to the $H^1(\Omega_f)-$norm i.e.
\begin{equation}\label{V}
\mathcal{V} := \overline{\mathcal{D}_f}^{H^1}, 
\end{equation}
the closed subspace $\mathcal{V}$ is endowed with the $H^1$-norm and thus it is a Hilbert space.

\subsection*{Well-posedness and energy estimate.}
\begin{proposition}\label{propoistion2}
$\forall F\in L^2(0,T; H^{-1}(\Omega_f))$, $\forall v_f^0\in L^2(\Omega_s)$ and  $\forall g_f\in L^2(0,T; H^{-\frac{1}{2}}(\Sigma))$, there  exists a unique solution $v_f$ to problem (\ref{sys4}) such that:
\begin{equation}\label{regulariteStokes}
v_f \in  L^2\left(0,T; \mathcal{V}\right)\cap H^1(0,T; H^{-1}(\Omega_f));
\end{equation}
moreover, one has the following energy estimate:
\begin{equation}\label{Ineq3}
\begin{aligned}
&  ||\partial_t v_f||_{L^2(0,T; H^{-1}(\Omega_f))} + ||v_f||_{L^2\left(0,T; \mathcal{V}\right)}\\& \leq C_f( ||\left( S(v_f, \pi) + a(t, \nabla v_f)\right)\cdot\overrightarrow{n} ||_{L^2\left( 0,T; H^{-\frac{1}{2}}(\Sigma)\right)}\\& + || F||_{L^2(0,T; H^{-1}(\Omega_f))}  + ||v_f^0||_{L^2(\Omega_s)}).
\end{aligned}
\end{equation}
\end{proposition} 
\noindent The well-posedness result stated in Proposition \ref{propoistion2} is rather classic. An equivalent weak formulation of (\ref{sys4}) can be derived by integrating the first equation in (\ref{sys4}), against $\phi\in L^2\left(0,T; \mathcal{V}\right)$, to obtain, cf. {\cite[Problem 3.15, p.371]{11}}:
\begin{equation}\label{eq5}
 (\partial_t v_f, \phi) +\mathcal{A}(v_f)(\phi)  - < (a(t, \nabla v_f) + S(v_f, \pi))\cdot\overrightarrow{n}, \phi>_{[\frac{1}{2}]^\ast, [\frac{1}{2}], \Sigma}  = <F, \phi   >,
\end{equation}
for a.e. time $0\leq t\leq T$, where
\begin{equation}\label{Af}
\mathcal{A}(v_f)(\phi) := \int_{\Omega_f} a(t, \nabla v_f)\nabla\phi{\rm d}x + \int_{\Omega_f}\varepsilon(v_f)\varepsilon(\phi){\rm d}x,
\end{equation}
where we used the fact ${\rm div}(\nabla v_f + \nabla^t v_f)=0 $ which holds since  ${\rm div}v_f = 0$, consequently one has ${\rm div}S(v_f, \pi) = \Delta v_f - \nabla\pi $ in this case, cf. {\cite[Problem 1.1, p.237-240]{6}}. Given the assumption (\ref{monotone}), the operator $\mathcal{A}$ is strongly monotone. Problem (\ref{eq5}) rewrites:

 \vspace{0.25cm}
 
 Find  $v_f \in L^2\left(0,T; \mathcal{V}\right)\cap H^1\left(0,T; H^{-1}(\Omega_f)\right)$ s.t.:
\begin{equation}\label{weakproblem}
\begin{aligned}
&-\int_0^T ( v_f, \partial_t\phi){\rm d}t + \int_0^T\mathcal{A}(v_f)(t)(\phi)(t) {\rm d}t \\&= \int_0^T <(a(t, \nabla v_f) + S(v_f, \pi))\cdot\overrightarrow{n}, \phi(t)>_{-\frac{1}{2}, \frac{1}{2}, \Sigma}{\rm d}t + \int_0^T  <F, \phi   >{\rm d}t, 
\end{aligned}
\end{equation}
$ \forall \phi\in L^2\left(0,T; \mathcal{V}\right)$ with $\partial_t\phi\in L^2\left(0,T; H^{-1}(\Omega_f)\right)$ and $\phi(T) =0$. Problem (\ref{weakproblem}) fits in the class of quasi-linear parabolic problems. One deals with such a problem using classical arguments, see for instance 
{\cite[Porposition 5.1, p.129]{5}} and  {\cite[Example 6.A, p.139]{5}}. Applying these last results we infer that problem (\ref{weakproblem}) admits a unique solution $v_f\in L^2\left(0,T; \mathcal{V}\right)$ for every $g_f\in L^2\left( 0,T; H^{-\frac{1}{2}}(\Sigma)\right)$. Regarding the existence issue, also cf.  {\cite[Theorem 1.1, p.225]{7}}. Furthermore, one can derive the energy estimate (\ref{Ineq3}) from problem (\ref{weakproblem}) by choosing as test function $\phi = v_f\in L^2(0,T; \mathcal{V}) $ and using the assumption (\ref{monotone}).

\begin{remark}\label{remark1}
 Let $v_f\in L^2\left(0,T; H^{\frac{1}{2}}(\Sigma)\right)\underset{continuous}{\hookrightarrow } L^2\left(0,T; L^2(\Sigma)\right) $. It is a classical fact that $v_f$ can be arbitrarily approximated by an element  $v_f^n \in C^\infty_0\left(0,T; H^{\frac{1}{2}}(\Sigma)\right)$ w.r.t. the norm of the space $L^2\left(0,T; H^{\frac{1}{2}}(\Sigma)\right)$. On another hand, by applying the Poincar\'e inequality in the time variable with $v_n$, we can easily show that:
\begin{equation}\label{p1}
 \int_0^t||u_f^n(r)||^2_{L^2(\Sigma)}{\rm d}r \leq C_p \int_0^t||v_f^n(r)||^2_{L^2(\Sigma)}{\rm d}r + 
 ||u_f^n(0)||_{L^2(\Sigma)}.    
\end{equation}
 On the other hand:
\begin{equation}\label{p2}
\begin{aligned}
&\int_0^T|u^n_f(t)|^2_{H^{\frac{1}{2}}(\Sigma)}{\rm d}t := \int_\Sigma\int_\Sigma\int_0^T \frac{ |u_f^n(t, x) - u_f^n(t, y)|^2}{|x - y|^2 }{\rm d}t {\rm d}x{\rm d}y\\& \leq C_p\int_\Sigma\int_\Sigma\int_0^T \frac{ |v_f^n(t, x) - v_f^n(t, y)|^2}{|x - y|^2 }{\rm d}t {\rm d}x{\rm d}y + T|u^n_f(0)|^2_{H^{\frac{1}{2}}(\Sigma)}.
\end{aligned}
\end{equation}
Combining estimates (\ref{p1}), (\ref{p2}) and using $u_f(0,x)=0 $ for $x\in\Sigma$, then letting $n\rightarrow\infty$: 
\begin{equation}\label{poincar}
\int_0^T||u_f(r)||^2_{H^{\frac{1}{2}}(\Sigma)}{\rm d}r \leq C_p\int_0^T||v_f(r)||^2_{H^{\frac{1}{2}}(\Sigma)}{\rm d}r.
\end{equation}
The same conclusion holds in case of fractional Sobolev spaces $H^{s}(\Sigma)$ with $s\in \mathbb{R}_+^*$.
\end{remark}
We infer from Remark \ref{remark1}, that the fluid displacement $u_f$ satisfies:
\begin{equation}\label{disreg}
u_f \in L^2\left(0,T;  H^\frac{1}{2}(\Sigma)\right)\cap H^1\left(0,T;  L^2(\Sigma)\right)\subset H^\frac{1}{2}((0,T)\times\Sigma).
\end{equation}
\subsection*{The idea.} Let us explain the main idea of this section. Define $T_2$ to be the operator: 
\begin{equation}\label{T2}
\begin{aligned}
T_2:   L^2\left( 0,T: H^{-\frac{1}{2}}(\Sigma)\right) &\rightarrow \mathcal{X}\subset H^\frac{1}{2}((0,T)\times\Sigma)\\
                 g_f &\mapsto  T_2(g_f) = u_f, 
\end{aligned}
\end{equation}
this operator associates to every Neumann data, on the fluid part of the contact interface, the displacement $u_f$ corresponding to the velocity $v_f$ which is a solution of (\ref{sys4}). It is easily seen, by combing (\ref{Ineq3}), (\ref{poincar}) and by applying (\ref{disreg}), that:
\begin{equation}\label{estimmixd}
\begin{aligned}
&||T_2(g_f^1) - T_2(g_f^2)||_{H^\frac{1}{2}\left((0,T)\times\Sigma\right)}\\&  =
||u_f^1 - u_f^2||_{H^\frac{1}{2}(0,T; L^2(\Sigma))}+ ||u_f^1 - u_f^2||_{L^2(0,T; H^\frac{1}{2}(\Sigma))} \\& \leq C_p  ||v_f^1 - v_f^2||_{L^2(0,T; H^1(\Omega_f))}\\&    \leq C_f ||g_f^1 - g_f^2||_{L^2\left( 0,T; H^{-\frac{1}{2}}(\Sigma)\right)}
\end{aligned}
\end{equation}
for every $g_f^1, g_f^2 \in L^2\left( 0,T; H^{-\frac{1}{2}}(\Sigma)\right)$, and thus the operator $T_2$  is continuous. Let $T_1$ and $T_2$ be defined respectively by (\ref{T1}) and (\ref{T2}). Define the operator:
\begin{equation}\label{T}
\begin{aligned}
T:   \mathcal{X} &\rightarrow \mathcal{X}\subset H^\frac{1}{2}((0,T)\times\Sigma)\\
                 u_s &\mapsto  T(u_s) := T_2\circ T_1(u_s) = u_f, 
\end{aligned}
\end{equation} 
We remark that the global solution of the coupled problem (\ref{sys1eq1})-(\ref{sys1eq7}) is a fixed point of $T$, then to show existence of a solution to (\ref{sys1eq1})-(\ref{sys1eq7}) it suffices to prove existence of a fixed point of the of the operator $T$. To be able to use fixed point theory we need some compactness. However, $T$ sends solid displacements from $:\mathcal{X}\subset H^\frac{1}{2}\left( (0,T)\times\Sigma\right)$ into no more spatial-regular fluid displacements, that is into: $\mathcal{X}\subset H^\frac{1}{2}\left( (0,T)\times\Sigma\right)$.  In order to recover some compactness we need to consider a sequence of auxiliary operators $T_2^\epsilon$. To achieve that, we proceed into a regularization of the Stokes problem (\ref{weakproblem}) i.e. to define a sequence of problems depending on a small real parameter, $\epsilon >0$, in such a way that the new operator $T_2^\epsilon\circ T_1$ sends $H^\frac{1}{2}\left((0,T)\times\Sigma\right)$ into a more regular space in the spatial variable, this will ensure the needed compactness. The (solution of the) original problem will be recovered by letting $\epsilon\rightarrow 0$.

\subsection*{Regularized problem.}\label{regz} Consider the space:
\begin{equation}\label{W}
 \mathcal{W} := \overline{\mathcal{D}_f}^{H^2}
\end{equation}
to be the completion of $\mathcal{D}_f$, defined by (\ref{Df}), with respect to the Sobolev $H^2-$norm. It is indeed a Hilbert spaces. We denote by $(. , .)_{H^2}$ its canonical inner product, and by $(.,.)_{H^2_{sn}}$ the part of $(. , .)_{H^2}$ that involves only the second derivatives. We denote $\mathcal{W}^\ast$ its dual. Given  $g_f \in L^2\left(0,T; H^{-\frac{1}{2}}  (\Sigma)\right)$, consider the regularized problem:

\vspace{0.25 cm}

Find   $v^\epsilon_f  \in  L^2\left( 0,T; \mathcal{W}\right)$ such that:
\begin{equation}\label{regularizedProblem}
\begin{aligned}
&\int_0^T ( \partial_tv_f^\epsilon, \phi)_{L^2}{\rm d}t  + \mathcal{A}_\epsilon(v_f^\epsilon)(\phi)    \\& = \int_0^T<F(t), \phi(t)>_{-1, 1, \Omega_f} {\rm d}t \int_0^T<g_f(t), \phi(t)>_{-\frac{1}{2}, \frac{1}{2}, \Sigma} {\rm d}t,
\end{aligned}
\end{equation}
for all $\phi \in L^2\left( 0,T; \mathcal{W} \right)$, where 
$$\mathcal{A}_\epsilon(v_f^\epsilon)(\phi) := \int_0^T \mathcal{A}(v_f^\epsilon)(\phi){\rm d}t +  \epsilon \int_0^T\int_{\Omega_f}(v^\epsilon_f , \phi)_{H^2_{sn}}~{\rm d}x{\rm d}t,$$
and where $\mathcal{A}$ is given by (\ref{Af}). It is easily seen that the quasi-linear elliptic operator $\mathcal{A}_\epsilon: \mathcal{W}\rightarrow \mathcal{W}^\ast$ satisfy the assumptions of {\cite[Proposition 5.1, p.129]{5}} and that of  {\cite[Theorem 5.1, p.128]{5}}. Then for every $g_f\in L^2\left(0,T; H^{-\frac{1}{2}}(\Sigma)\right)$, problem (\ref{regularizedProblem}) admits a unique solution:
\begin{equation}\label{moreref}
v_f^\epsilon \in L^2\left(0,T; \mathcal{W}\right)
\end{equation}
for all $\epsilon>0$, where $\mathcal{W}$ is given by (\ref{W}). An energy estimate can be derived by using the strong monotony of the operator $\mathcal{A}_\epsilon$. We thus obtain the estimate:
\begin{equation}\label{estimateStokesreg0}
||v_f^{1, \epsilon} - v_f^{2,\epsilon}||_{L^2\left(0,T; \mathcal{W}\right)}\leq C_f ||g^1_f - g^2_f ||_{L^2\left(0,T; H^{-\frac{1}{2}}(\Sigma) \right)}.
\end{equation}
Estimate (\ref{estimateStokesreg0}) is obtained by choosing as test function $\phi = v_f^{1,\epsilon} - v_f^{2, \epsilon}$, in (\ref{regularizedProblem}), corresponding to $g_f^1, g_f^2$.
 Using the continuity of the trace operator $$\gamma_0: H^2(\Omega_f)\rightarrow H^{\frac{3}{2}}(\Sigma),$$ we infer from (\ref{estimateStokesreg0}):
 \begin{equation}\label{estimateStokesreg1}
 \begin{aligned}
&||v_f^{1,\epsilon} - v_f^{2,\epsilon}||_{L^2\left(0,T; H^{\frac{3}{2}}(\Sigma)\right)}\\&\leq C_f ||g^1_f - g^2_f||_{L^2\left(0,T; H^{-\frac{1}{2}}(\Sigma) \right)}
\end{aligned}
\end{equation}

\begin{remark}\label{remark2}
The same fact that was stated in Remark \ref{remark1} holds with $u_f^\epsilon$ and $v_f^\epsilon$ in the space $L^2\left(0,T; H^{\frac{3}{2}}(\Sigma)\right)$. Thus we obtain by combining (\ref{poincar}) corresponding to the space $H^{\frac{3}{2}}(\Sigma)$ and estimate (\ref{estimateStokesreg1}):
 \begin{equation}\label{Ineq4}
 \begin{aligned}
&||u_f^{\epsilon,1} - u_f^{\epsilon,2}||_{L^2\left(0,T; H^1(\Sigma)\right)} + ||u_f^{\epsilon,1} - u_f^{\epsilon,2}||_{H^1\left(0,T; L^2(\Sigma)\right)} \\&\leq C_f ||g^1_f - g^2_f ||_{L^2\left(0,T; H^{-\frac{1}{2}}(\Sigma) \right)},
\end{aligned}
\end{equation}
for every $g_f^1, g_f^2\in   L^2\left( 0,T: H^{-\frac{1}{2}}(\Sigma)\right).$ 
 \end{remark}
\subsection*{The auxiliary operator $T^\epsilon$.} First define $\mathcal{Z}$ to be the space:
\begin{equation}\label{Z}
\mathcal{Z} := H^1\left(0,T; L^2(\Sigma)\right)\cap  L^2\left(0,T; H^1_0(\Sigma)\right)\subset H^1((0,T)\times\Sigma),
\end{equation}
it is endowed with the standard $H^1((0,T)\times\Sigma)-$norm. Define the operator $T^\epsilon_2$ by: 
 \begin{equation}\label{T2epsilon}
\begin{aligned}
T^\epsilon_2:   L^2\left( 0,T; H^{-\frac{1}{2}}(\Sigma)\right) &\rightarrow \mathcal{Z}\subset\mathcal{X}\\
                 g_f &\mapsto  T_2^\epsilon(g_f) = u_f^\epsilon, 
\end{aligned}
\end{equation}
i.e. it associates to every $g_f \in L^2\left( 0,T; H^{-\frac{1}{2}}(\Sigma)\right)$, the trace on $(0,T)\times\Sigma$ of the fluid displacement $u^\epsilon_f$ corresponding to the fluid velocity $v_f^\epsilon$ which is the solution of the regularized problem (\ref{regularizedProblem}) with $g_f$ as a Neumann data. Estimate (\ref{Ineq4}) translates the boundedness of $T^\epsilon_2$. Let us introduce the auxiliary operator $T^\epsilon := T^\epsilon_2 \circ T_1$, that is: 
\begin{equation}\label{Tepsilon}
\begin{aligned}
T^\epsilon: \mathcal{X}&\rightarrow \mathcal{Z}\subset \mathcal{X}\\
                 u_s &\mapsto u_f^\epsilon = T_2^\epsilon\circ T_1(u_s),
\end{aligned}
\end{equation}
where $\mathcal{X}$ and $\mathcal{Z}$ are respectively given by (\ref{X}) and (\ref{Z}), moreover the operators $T_1$ and $T_2^\epsilon$ are respectively defined by (\ref{T1}) and (\ref{T2epsilon}). The operator $T^\epsilon$ associates to every solid displacement $u_s \in  \mathcal{X}$, the fluid displacement $u_f^\epsilon\in \mathcal{X}$ defined on the fluid part of the contact interface $(0,T)\times\Sigma$. 

The auxiliary coupled problem is formed out of the Lam\'e solid problem (\ref{sys2}) and the quasi-linear Stokes regularized problem (\ref{regularizedProblem}). By combining the estimates (\ref{estimlame}) and (\ref{Ineq4}) on one hand, and using the coupling condition (\ref{eq1}) on the other hand, we infer the estimate:
\begin{equation}\label{continuT}
 \forall\epsilon>0,\quad ||T^\epsilon(u_s^1) - T^\epsilon(u_s^2)||_{\mathcal{X}}\leq C_s C_f ||u_s^1 - u_s^2||_{\mathcal{X}}~,
\end{equation}
for all $u_s^1, u_s^2\in \mathcal{X}$ i.e. the non-linear operator $T^\epsilon$ is Lipschitz for every $\epsilon >0$. One should notice that the constant $ C_sC_f$ doesn't depend on $\epsilon$.

\section{Existence of solution for the main coupled problem } 
\subsection{ Existence result for the auxiliary coupled problem.}
\label{ss41}
Let $T^\epsilon$ be as defined by (\ref{Tepsilon}) and $\mathcal{X} $
be as defined by (\ref{X}).  We propose to prove the following proposition:
\begin{proposition}\label{propositon3}
For every $\epsilon>0$, the operator $T^\epsilon$ admits a fixed point $u^\epsilon\in \mathcal{X}$.
\end{proposition}
The idea for proving Proposition \ref{propositon3} is to combine the Banach and the Schaefer fixed point theorems. Let us recall Schaefer's theorem, cf. {\cite[Theorem 4.3.2, p.29]{8}}:
\begin{theorem}\label{Schaefer}(Schaefer). Let $\mathcal{K}$ be a Banach space. Let $T:\mathcal{K} \rightarrow \mathcal{K}$ be a continuous and compact mapping. Assume that the set
$$ \lbrace  v\in \mathcal{K}: ~v = \rho Tv~\text{for some}~\rho,~0 \leq \rho\leq 1   \rbrace   $$
is bounded, then $T$ admits a fixed point in $\mathcal{K}$.
\end{theorem}

We are going now to prove Proposition \ref{propositon3}:

\noindent\begin{proof} We deal at first with a contraction mapping. Fix $\rho \in\mathbb{R}$ such that $\rho C_sC_f <1$. The coefficient $\rho$ doesn't depend on $\epsilon$. Consider the operator $\rho T^\epsilon$ defined by:
\begin{alignat*}{2}
\rho T^\epsilon:  \mathcal{X} &\rightarrow  \mathcal{X}\cap\mathcal{Z}\subset \mathcal{X}\\
                 u_s &\mapsto u_f^\epsilon =\rho T_2^\epsilon\circ T_1(u_s),
\end{alignat*}
where $\mathcal{Z}$ is defined by (\ref{Z}). Applying estimate (\ref{continuT}), we have: 
\begin{equation}\label{ineq5}
 ||\rho T^\epsilon(u_s^1) - \rho T^\epsilon(u_s^2)||_{\mathcal{X}} \leq \rho C_sC_f ||u_s^1 - u_s^2||_{\mathcal{X}},
\end{equation}
for all $u_s^1, u_s^2\in \mathcal{X}$. Let us note the following four facts:
\begin{itemize}
\item {\it  Continuity of the operators $T^\epsilon$.} Estimate (\ref{continuT}) means that $T^\epsilon$ is Lipschitz and thus continuous. Moreover, estimate (\ref{ineq5}) implies that $\rho T^\epsilon$ is a contraction mapping for every $\epsilon>0$. Consequently, since $\mathcal{X}$ is a Banach space, we deduce, by applying the Banach fixed point theorem, that $\rho T^\epsilon$ admits a unique $[0,T]-$globally defined fixed point $u^\epsilon_\rho \in \mathcal{X}$.
 
 \item {\it Boundedness condition}. We deduce from the preceding point that, $\forall\epsilon>0$, the set $ \lbrace  u\in \mathcal{X}, ~u = \rho T^\epsilon u~  \rbrace   $ reduces to the unique fixed point $u^\epsilon_\rho$, and thus it is bounded.
 
\item  {\it  Stability of the operators $ T^\epsilon$.} According to definition (\ref{Tepsilon}), we have $T^\epsilon \left( \mathcal{X} \right) \subset \mathcal{X}$, thus $T^\epsilon$ is stable for every $\epsilon>0$.

\item {\it  Compactness of the operator $ T^\epsilon$.} Let $\epsilon>0$. By combining estimates (\ref{estimlame}), (\ref{Ineq4}) and the boundedness of $||T^\epsilon(0)||_{\mathcal{Z}}$, we remark that the operator $T^\epsilon$ sends every bounded subset $E\subset \mathcal{X}$ into an $H^1-$bounded subset $T^\epsilon(E)\subset H^1((0,T)\times\Sigma)$. Using the compact embedding $H^1((0,T)\times\Sigma)  \underset{compact}{\hookrightarrow } H^\frac{1}{2}((0,T)\times\Sigma)$, we deduce that $T^\epsilon(E)$ is compact in $\mathcal{X}$, thus $T^\epsilon$ is compact $\forall\epsilon >0$.
\end{itemize}

\noindent Thus we checked for $T^\epsilon$ the sufficient conditions of Theorem \ref{Schaefer}; this proves Proposition \ref{propositon3}.
\end{proof}

\subsection{Existence of a solution to the coupled problem (\ref{sys1eq1})-(\ref{sys1eq7}).} We now establish the existence of a solution for the global coupled problem which, as pointed out above, amounts at establishing the existence of a fixed point of the operator $T$ given by (\ref{T}). According to Proposition \ref{propositon3}, $\forall\epsilon >0$, there exists $u^\epsilon\in\mathcal{X}$ such that 
\begin{equation}\label{fixedpointT}
T^\epsilon u^\epsilon = u^\epsilon,
\end{equation} 
where $\mathcal{X}$ is given by (\ref{X}). Since $(u^\epsilon)_{\epsilon>0}$ is an uncountable family of vectors belonging to the separable normed vector space $(\mathcal{X}, ||~||_{\mathcal{X}})$, then it must have a limit point with respect to the topology induced by the norm of $\mathcal{X}$ i.e. there exists $u^0\in\mathcal{X}$ such that:   
\begin{equation}\label{convergence1}
||u^\epsilon - u^0||_{\mathcal{X}}\rightarrow 0, 
\end{equation}
as $\epsilon\rightarrow 0$. To be able to pass to the limit $\epsilon\to 0 $ in (\ref{fixedpointT}), we need Lemma \ref{lemma1}:
\begin{lemma}\label{lemma1}
Let $T^\epsilon$ and $T$ be respectively as defined  by (\ref{Tepsilon}) and (\ref{T}), then:
\begin{equation}\label{converegnceTepsilon}
 T^\epsilon u    \underset{strongly~in~\mathcal{X}}{\rightarrow}        Tu ,\quad  \forall u\in \mathcal{X}.
\end{equation}
\end{lemma}

\begin{proof} 
 Fix $u\in\mathcal{X}$ and pose $g_f = T_1 u$, where $T_1$ is given by (\ref{T1}). Choose in (\ref{regularizedProblem}) $\phi = u_f^\epsilon$ where $u_f^\epsilon = T_2^\epsilon(g_f)$ and where $T_2^\epsilon$ is given by (\ref{T2epsilon}). Using (\ref{Ineq4}), we infer immediately that  $T^\epsilon u    \underset{weakly~in~\mathcal{X}}{\rightharpoonup}   U_f$ with $U_f\in \mathcal{X}$, let us show that $U_f = Tu $.
  
Denote $v_f^\epsilon$ to be the velocity  $v_f^\epsilon := \partial_t u_f^\epsilon$. It is clear, by using (\ref{estimateStokesreg0}) with $g_f^2 = 0$ and $g_f^1 = g_f$, that $\left(||v_f^\epsilon||_{L^2(0,T; H^1(\Omega_f))}\right)_\epsilon $ is uniformly bounded, thus: 
\begin{equation}\label{weaklimit1}
 v^\epsilon_f \underset{weakly}{\rightharpoonup}V_f \quad{\rm in}~L^2(0,T; H^1(\Omega_f)).
\end{equation}
Moreover, by combining (\ref{estimateStokesreg0}) and (\ref{weaklimit1}) we infer that both $\epsilon\int_0^T |v_f^\epsilon(t)|^2_{H^2_{sn}(\Omega_f)}{\rm d}t:= \epsilon\int_0^T(v_f^\epsilon, v_f^\epsilon)_{H^2_{sn}}{\rm d}t$ and $||\partial_t v_f^\epsilon||_{L^2(0,T; H^{-1}(\Omega_f))}$
 are also uniformly bounded with respect to $\epsilon>0$, consequently:
\begin{equation}\label{weaklimit2}
\begin{aligned}
\partial_t v^\epsilon_f &\underset{weakly}{\rightharpoonup}\partial_t V_f  \quad{\rm in}~L^2(0,T; H^{-1}(\Omega_f)), \\
\epsilon\partial^2_{x_ix_j} v^\epsilon_f &\underset{weakly}{\rightharpoonup} 0 \quad{\rm in}~L^2(0,T; L^2(\Omega_f)).
\end{aligned}
\end{equation}

Combining the weak problems (\ref{weakproblem}), (\ref{regularizedProblem}) and applying the last convergence in (\ref{weaklimit2}) on one hand, and by considering the auxiliary function $u^n$ corresponding to problem (\ref{weakproblem}) with $f^n$ an d$g^n$ as  right hand sides such that $||f^n - f||\to 0$ and $ ||g^n - g||\to 0$, we deduce that:
\begin{equation}\label{strongconv}
||\nabla v_f^\epsilon - \nabla v_f||_{L^2(0,T; L^2(\Omega_f))} \to 0, 
\end{equation}
it follows that $V_f = v_f$. Given that $a$ is Lipschitz in the second variable, then (\ref{strongconv}) yields:
\begin{equation}\label{consa}
  a(t, \nabla v_f^\epsilon ) \to a(t, \nabla V_f)  \quad{\rm converges ~~ strongly}.
\end{equation}
Finally, letting $\epsilon\rightarrow 0$ in (\ref{regularizedProblem}) by mean of combining (\ref{weaklimit1}), (\ref{weaklimit2}) and (\ref{consa}), one infers that $V_f$ is the fluid velocity field corresponding to the solid displacement $u\in\mathcal{X}$. Consequently we have $\partial_t u  = V_f  $ and $U_f = Tu$, this yields convergence (\ref{converegnceTepsilon}).
\end{proof}

\subsection*{Proof of Theorem \ref{main-The}}
Now we are ready to present a proof of the main result:
\noindent\begin{proof} 
Let $u^0$ be such as defined by (\ref{convergence1}), we have: 
\begin{equation}\label{lim1}
\begin{aligned}
& || T^\epsilon u^\epsilon - T u^0||_{\mathcal{X}} \\&  \leq || T^\epsilon u^\epsilon - T^\epsilon u^0||  +  || T^\epsilon u^0  - T u^0||.
\end{aligned}
\end{equation}

On one hand we have:
\begin{equation*}
\begin{aligned}
  &  || T^\epsilon u^\epsilon - T^\epsilon u^0||_{\mathcal{X}} \\& \leq  C_sC_f||u^\epsilon - u^0||_{\mathcal{X}},  
  \end{aligned}
\end{equation*}
where we used (\ref{continuT}). Applying (\ref{convergence1}), we infer:
 \begin{equation}\label{lim3}
\begin{aligned}
  &  || T^\epsilon u^\epsilon - T^\epsilon u^0||_{\mathcal{X}}\to 0\quad{\rm as}~\epsilon\to 0.
  \end{aligned}
\end{equation} 

On another hand, using (\ref{converegnceTepsilon}), we obtain:
\begin{equation}\label{lim4}
\begin{aligned}
  &  || T^\epsilon u^0  - T^0 u^0||\rightarrow 0,~ \epsilon\rightarrow 0.
\end{aligned}
\end{equation}

 Combining (\ref{lim3}) and (\ref{lim4}) it yields 
 \begin{equation}\label{limit22}
 || T^\epsilon u^\epsilon - T u^0||_{\mathcal{X}},~ \epsilon\rightarrow 0.
 \end{equation}
 
Combining (\ref{fixedpointT}), (\ref{convergence1}) and (\ref{limit22}), we deduce:
\begin{equation}\label{eq6}
u^0 = \lim_{\epsilon\rightarrow 0} u^\epsilon   = \lim_{\epsilon\rightarrow 0} T^\epsilon u^\epsilon  =  T u^0,
\end{equation}
that is $T$ has a fixed point $u^0 \in \mathcal{X}$. This completes the proof of Theorem \ref{main-The}. 

Regarding the regularity claimed in Theorem \ref{main-The}, one can use estimate (\ref{Ineq1}) to infer the regularity of the solid displacement $u_s\in\mathcal{X}$ on $\Omega_s$. Furthermore, one combines the estimates (\ref{regulariteStokes}) and (\ref{estimlame}) to infer the regularity of the fluid velocity $v_f$ on $\Omega_f$. The regularity of the pressure $\pi$ can be inferred in the following fashion: one considers the regularity of the solution $u^0$ stated in (\ref{regulariteStokes}), namely that  $\partial_t u^0 \in L^2\left(0,T; H^{-1}(\Omega_f)\right)$ and then applies the energy estimate in {\cite[Theorem 25, p.226]{9}} and thus infer that $\pi\in L^2\left( 0,T; L^2(\Omega_f)/\mathbb{R}\right)$.
\end{proof}

\section{Conclusion}
\noindent According to the procedure adopted by the authors in \cite{1}, cf. the bottom of {\cite[Proof of Theorem 5.1, p.571]{1}}, the globally-in-time existence of a solution can be inferred by using an iterative method which is based on the linearity of the problem they considered. We emphasize that this method is no longer applicable in the non-linear context of the present work.

It should be noted that there is no memory effect in the (iterative) resolution of the presently addressed coupled problem in the sense that the result (in particular the displacement) obtained at a time $ t $ does not depend on the displacement at an earlier time. This is notably due to the fact that the deformation of the geometry (in particular that of the contact surface $ \Sigma $), within time incrementing i.e. during the coupling, is not taken into account, although we assumed in the introduction that the domain $\Omega_f(t)$ depends upon time.

 That said, the main theorem established in this paper remains relevant. Indeed, it can be incorporated as an auxiliary result to demonstrate a more complex well-posedness result like for instance {\cite[Theorem 2.1, p.555]{1}}.  If one is willing to prove a result analogous to the later in the framework considered in this paper, then one has, in a first step, to generalize Theorem \ref{main-The} to the case when the data in system (\ref{sys1eq1})-(\ref{sys1eq7}) are non-homogeneous along with a non necessary free divergence condition. We believe that this step can be achieved with some slight modification of the procedure described in the present work, one obtains up to this step a result analogous to {\cite[Theorem 5.1, p.570]{1}. Then one should follow the same method as in the proof of  {\cite[Theorem 2.1, p.555]{1}} which deals with the well-posedness of system {\cite[Problem 2.1, p.551-552]{1}} and which requires taking into account the deformation of the geometry during the interaction, especially that of the contact interface.

\end{document}